\documentclass[12pt,leqno,a4paper]{article}
\usepackage{mathrsfs}

\usepackage{amsfonts}
\usepackage{amsmath, amssymb, amsthm, amsfonts}

\numberwithin{equation}{section}

\def\txt#1{{\textstyle{#1}}}
\def\hf{{\textstyle{\frac12}}}
\def\a{\alpha}\def\b{\beta}
\def\d{{\,\rm d}}
\def\e{\varepsilon}

\def\G{\Gamma} 

\def\s{\sigma}
\def\z{\zeta}
\def\t{\theta}
\def\={\;=\;}
\def\le{\leqslant}
\def\ge{\geqslant}

\def\zt{\zeta(\hf+it)}

\begin{document}
\baselineskip=17pt
\title{\large \bf Hardy's function $Z(t)$ - results and problems}
\author{\bf Aleksandar Ivi\'c$^1$}
\date{}
\bigskip

\footnotetext[0]{$^1$Serbian Academy of Science and Arts, Knez Mihailova 35,
11000 Beograd, Serbia\\
e-mail: {\tt aleksandar.ivic@rgf.bg.ac.rs, aivic\_2000@yahoo.com}}
%\\

 \maketitle

{\small {\bf Abstract}.
This is primarily an overview article
on some results and problems involving the classical Hardy function
$$
Z(t) := \zt{\bigl(\chi(\hf+it)\bigr)}^{-1/2}, \quad \z(s) = \chi(s)\z(1-s).
$$
In particular, we discuss the first and third moment of $Z(t)$ (with and without
shifts) and the
distribution of its positive and negative values. A new result involving the
distribution of its values is presented.}


\bigskip

\section{\bf Definition of Hardy's function}

%\subsection{The definitions}

The primary aim of this paper is to present some results and problems involving
{\it Hardy's function} $Z(t)$, since in recent years there was a revival of
interest in its study.
This classical  function (see e.g., the author's monograph \cite{Iv7} for an extensive account)
has a century long history. It is defined as
$$
Z(t) := \zt{\bigl(\chi(\hf+it)\bigr)}^{-1/2},\leqno(1.1)
$$
where $\chi(s)$ comes from the familiar functional equation for $\z(s)$
(see e.g., Chapter 1 of \cite{Iv1}),
namely $\z(s) = \chi(s)\z(1-s)$ for $s\in \mathbb C$, so that
$$
\chi(s) = 2^s\pi^{s-1}\sin(\hf \pi s)\G(1-s),\quad
\chi(s)\chi(1-s)=1.
$$
It follows  that
$$
\overline{\chi(\hf + it)} = \chi(\hf-it)= \chi^{-1}(\hf+it),
$$
so that $Z(t)\in\mathbb R$ when $t\in\mathbb R$, $Z(t) = Z(-t)$,
and $|Z(t)| =|\zt|$. Thus the zeros of $\z(s)$ on the ``critical line'' $\Re s =1/2$
correspond to the real zeros of $Z(t)$, which makes $Z(t)$ an invaluable tool
in the study of the zeros of the zeta-function on the critical line.
Alternatively, if we use the symmetric form of the functional equation
for $\z(s)$, namely
$$
\pi^{-s/2}\z(s)\G(\hf s) = \pi^{-(1-s)/2}\z(1-s)\G(\hf(1-s)),
$$
then for $t\in \mathbb R$ we obtain
$$
Z(t) = {\rm e}^{i\t(t)}\zt,\quad {\rm e}^{i\t(t)} := \pi^{-it/2}\frac
{\G(\frac{1}{4}+\hf it)}{|\G(\frac{1}{4}+\hf it)|}\quad (\t(t) \in\mathbb R),
$$
which implies that $Z(t)$ is a smooth function.

\medskip
For completeness, recall that the Riemann zeta-function is defined by
$$
\z(s) = \sum_{n=1}^\infty n^{-s} = \prod_{p}{(1-p^{-s})}^{-1}\leqno(1.2)
$$
for $\Re  s > 1$, where $p$ denotes primes.
For other values of the complex variable $s = \s+it\, (\s,t \in\mathbb R)$ it is defined
by  analytic continuation.
It is regular for $s\in\mathbb C$, except at $s=1$ where it has a simple
pole with residue 1.
The product representation in (1.2) shows that $\z(s)$ does not vanish for $\s>1$.
The best known ``zero-free region'' for $\z(s)\;$ is of the form
$$
\s > 1 - C(\log t)^{-2/3}(\log\log t)^{-1/3}\quad(C>0, t \ge t_0 >0).\leqno(1.3)
$$
This was obtained in 1958 by the method of I.M. Vinogradov
(see e.g., his  works \cite{Vin}, \cite{Vin1}, Chapter 4 of \cite{KaVo} and Chapter 6 of \cite{Iv1}).
The best known numerical
values in (1.3) are $C = 1/57,54,\, t_0 =4$, and they are due to K. Ford \cite{For}.

\section{Zeta-zeros on the critical line }
Hardy's original application \cite{Hard} in 1914 of $Z(t)$
was to show that $\z(s)$ has infinitely
many zeros on the critical line $\Re s =1/2$  (see e.g., E.C. Titchmarsh [23]).
The argument is briefly as follows.
Suppose on the contrary that, for $T\ge T_0$, the function $Z(t)$ does not change sign. Then
$$
\int_T^{2T}|Z(t)|\d t \;=\; \Bigl|\int_T^{2T}Z(t)\d t\Bigr|.\leqno(2.1)
$$
On one hand we have
$$
\int_T^{2T}|Z(t)|\d t =\int_T^{2T}|\zt|\d t \ge \left|\int_T^{2T}\zt\d t\right|.
$$
One has the elementary formula (see e.g., Chapter 1 of \cite{Iv8})
\[
\zt = \sum_{n\le T}n^{-1/2-it} + \frac{T^{1/2-it}}{it-1/2} +O(T^{-1/2})\qquad(T\le t \le2T).
\leqno(2.2)
\]
Using (2.2) it is easily found, on integrating termwise the right-hand side, that
$$
\int_T^{2T}\zt\d t = 2T + O(T^{1/2}).
$$
This yields
$$
\int_T^{2T}|Z(t)|\d t \;=\;\int_T^{2T}|\zt|\d t \;\gg\;T,\leqno(2.3)
$$
and the even slightly sharper lower bound (see K. Ramachandra \cite{Ram}) $T(\log T)^{1/4}$ holds.

On the other hand, to bound the integral on the right-hand side of (2.1)
we can use the approximate functional equation (this is a weakened
form of the so-called Riemann--Siegel formula; for a proof see \cite{Iv1}  or \cite{Tit})
$$
Z(t) = 2\sum_{n\le\sqrt{t/(2\pi)}}n^{-1/2}\cos\Biggl(t\log\sqrt{\frac{t/(2\pi)}{n}}
- \frac{t}{2} - \frac{\pi}{8}\Biggr) + O\Bigl(\frac{1}{t^{1/4}}\Bigr).\leqno(2.4)
$$
If this expression is integrated and the second derivative test is applied
(see \cite{Iv1}  or \cite{Tit})
it follows that
$$
\int_T^{2T}Z(t)\d t \;=\; O(T^{3/4}).\leqno(2.5)
$$
Thus from (2.1)--(2.5) we obtain
$$
T \;\ll\; \int_T^{2T}|Z(t)|\d t \;\ll\;T^{3/4},
$$
which is a contradiction. This proves that $\z(s)$ has infinitely
many zeros on the critical line. In fact, the argument that leads to (2.5) actually shows that
$$
N_0(T) \;\gg\; T^{1/4},
$$
where $N_0(T)$ denotes the number of complex zeros $\rho = \b +i\gamma$ of $\z(s)$ for which
$\b = 1/2, 0<\gamma \le T$.  Later Hardy refined his argument to show that $N_0(T) \gg T$.

A. Selberg (see \cite{Sel} or \cite{Tit}
for a proof) improved this bound to $N_0(T) \gg T\log T$, which is one of the most important results
of analytic number theory. In fact, this implies that
$$
N_0(T) \;\ge\;CN(T) \leqno(2.6)
$$
for some $C>0$ and $T\ge T_0>0$. Here $N(T)$ denotes the number of $\rho = \b +i\gamma$
for which $0<\gamma \le T$. One has the classical {\it Riemann--von Mangoldt formula}
(see e.g., Chapter 1 of \cite{Iv1} for a proof)
\[
N(T) = \frac{T}{2\pi}\log \frac{T}{2\pi} - \frac{T}{2\pi} + O(\log T),
\]
and therefore (2.6) holds as a consequence of $N_0(T) \gg T\log T$.

N. Levinson \cite{Lev} in 1974 showed that $C = 1/3$ is permissible  in (2.6), J.B. Conrey \cite{Con}
1989 obtained $C=2/5$, that is, 40\% of the zeta-zeros are on the critical line.
The latest record was achieved by S. Feng \cite{Fen},
who proved that at least 41.73\% of the zeros of $\z(s)$
are on the critical line and at least 40.75\% of  those zeros
are simple ($\z(\rho) =0 \Rightarrow \z'(\rho)\ne 0$) and on the critical line.
Selberg's proof of (2.5) involved combining a ``mollifier''
 to compensate for irregularities in the
size of $|\z(s)|$ and the method of Hardy (and Littlewood).
Levinson introduced new ideas,
and subsequent research refined on the existing methods.

\medskip
{\bf Notation}. Owing to the nature of this text, absolute consistency in notation could
not be attained, although whenever possible standard notation is used. By
$\mathbb N,\mathbb Z,\mathbb R,\mathbb C$ we denote the set of natural numbers, integers, real and complex
numbers, respectively. The symbol $\e$ will denote arbitrarily small positive numbers,
not necessarily the same ones at each occurrence. The Landau symbol $f(x) = O\bigl(g(x)\bigr)$
and the Vinogradov symbol $f(x) \ll g(x)$ both mean that $|f(x)| \le Cg(x)$
for some constant $C>0, g(x)>0$ and
$x \ge x_0 > 0$. By $f(x) \ll_{a,b,\ldots} g(x)$ we mean that the constant implied by
the $\ll$-symbol depends on $a,b,\ldots$. The symbol $f(x) = \Omega_\pm\bigl(g(x)\bigr)$ means that both
$\limsup\limits_{x\to\infty}f(x)/g(x) > 0$ and $\liminf\limits_{x\to\infty}f(x)/g(x) < 0\,$ holds.

\section{Moments of Hardy's function}
\subsection{Discussion of $F_k(T)$}
For $k \in \mathbb N$ fixed, consider the $k$-th moment of $Z(t)$, namely the integral
$$
F_k(T) := \int_0^T Z^k(t)\d t.\leqno(3.1)
$$
Since $|Z(t)| = |\zt|$, it transpires that
$$
F_{2k}(T) \equiv \int_0^T|\zt|^{2k}\d t,
$$
which is one of the fundamental objects in the study of $\z(s)$. Even moments in general
are a natural object of study, because of the elementary identity $|z|^2 = z\cdot {\bar z}$.
When $z = \zt^k$, this permits one to develop the square and use various approximate
functional equations etc.
The reader is referred to the monographs of K. Ramachandra \cite{Ram} and the author
\cite{Iv2}, which deal exclusively with mean values (moments) of $\z(s)$. Also the
books of E.C. Titchmarsh \cite{Tit} and the author \cite{Iv1} contain a lot of material
on this subject, as does his review paper \cite{Iv8}.
Thus only the study of $F_{2k-1}(T)$ represents a novelty. The function
$Z(t)$ takes positive and negative values (and, heuristically, with a certain regularity),
so that one expects there will be a lot of cancellations when one evaluates $F_{2k-1}(T)$.
However, the following natural problem seems challenging in the general case.

\medskip
{\bf Problem 1}. {\it Show that, for $k>1$ a fixed integer, one has}
$$
F_{2k-1}(T) \,=\,o\left(\int_0^T|\zt|^{2k-1}\d t\right)\qquad(T\to\infty).\leqno(3.2)
$$

\medskip
The use of the Riemann--Siegel formula (2.4) does not seem adequate in proving (3.2).
Although we have at our disposal {\it smooth variants} of this formula, which
will be discussed a little
later, the problem of establishing (3.2)  nevertheless remains open. For a discussion
involving problems with $Z(t)$, see the author's paper \cite{Iv5}.

\medskip
\subsection{Bounds for  $F_1(t)$}
We turn now to $F(T) \equiv F_1(T)$. In 2004 the author \cite{Iv3}
improved (2.5) by obtaining a much stronger result than (3.2) for $k=1$, namely

\medskip
{\bf Theorem 1}. {\it We have}
$$
F(T) = \int_0^T Z(t)\d t = O_\e(T^{1/4+\e}).\leqno(3.3)
$$

\medskip
We sketch briefly the proof of (3.3). It is based  (A.I. \cite{Iv2}, 1990) on the use
of a  smooth approximate functional equation for $Z^k(t)$, namely
$$
Z^k(t) = 2\sum_{n\le2\tau}\rho\left(\frac{n}{\tau}\right)d_k(n)n^{-1/2}
\cos\left(t\log\frac{\tau}{n} - \frac{k}{2}t - \frac{\pi k}{8}\right)
+ O(t^{\frac{k}{4}-1}\log^{k-1}t),
$$
where for any fixed integer $k\ge1$,  $t\ge 2$,
$$
\tau \;=\; \left(\frac{t}{\pi} \right)^{k/2},
$$
and further notation is as follows. The generalized divisor function $d_k(n)$
(generated by $\z^k(s)$, which makes it possible to define $d_k(n)$
for an arbitrary $k\in\mathbb C$) represents
the number of ways $n$ may be represented as the product of $k$
factors ($d_1(n) \equiv 1, d_2(n) \equiv d(n)$, the number of divisors
of $n$). The test function $\rho(x)$ is a non-negative, smooth function
supported in $\,[0,2]\,$, such that $\rho(x) = 1$ for $0 \le x\le 1/b$
for a fixed constant $b>1$, and
$\rho(x) + \rho(1/x) = 1$ for all $x$. The last condition induces a
symmetry in the approximate functional equation for $Z^k(t)$.

For $k=1$ the error term gives $O(T^{1/4})$ after integration. The integration
of the main term
produces exponential integrals which are evaluated by the classical saddle point
method. There are a number of such results in the literature (see e.g., Chapter 2 of \cite{Iv1}).
The one that is convenient is the following lemma
(see p. 71 of the monograph by Karatsuba--Voronin \cite{KaVo}).

\medskip
{\bf Lemma 1}.
{\it If $f(x) \in C^{(4)}[a,b], f''(x) > 0$ in $[a,b]$, then
\begin{align*}&
\int_a^b e^{2\pi if(x)}\d x =
e^{\pi i/4}\frac{e^{2\pi if(c)}}{\sqrt{f''(c)}} \\&+\, O(AV^{-1})
+\,O\Bigl(\min({|f'(a)|}^{-1},\sqrt{A}\,)\Bigr) + 
O\left(\min({|f'(b)|}^{-1},\sqrt{A}\,)\right),
\end{align*}
where the main term is to be halved if $c=a$ or $c=b$, and
\begin{align*}
&0 < b - a \le V,\; f'(c) =0, \;a\le c \le b,\\&
f''(x) \asymp A^{-1},\; f^{(3)}(x) \ll (AV)^{-1},\; f^{(4)}(x) \ll A^{-1}V^{-2}\quad(A>0).
\end{align*}
}

Application of Lemma 1 and subsequent estimations and simplifications lead
 eventually to the upper bound in (3.3).

\medskip
In \cite{Iv3} it was conjectured that
$$
\int_0^T Z(t)\d t = O(T^{1/4}), \quad \int_0^T Z(t)\d t = \Omega_\pm(T^{1/4}).\leqno(3.4)
$$
This was proved, independently and by different methods,
 by M. Jutila \cite{Jut1}, \cite{Jut2} and M.A. Korolev \cite{Kor}.
Therefore they established (up to the value of the numerical constants which are involved
in the $O$ and $\Omega_\pm$ symbols) the true
order of the integral in question. For the integral in (3.4)
Korolev actually obtained the explicit bound
\[
\left|\int_{2\pi}^T Z(t)\d t\right| < 18.2T^{1/4}\qquad(T\ge T_0).
\]
\medskip
\subsection{The cubic moment of $Z(t)$}

In what concerns the cubic moment $F_3(T) \equiv \int\limits_0^TZ^3(t)\d t$, in Oberwolfach 2003 I
posed the following

\medskip
{\bf Problem 2}. {\it Does there exist a constant $0<c<1$
such that
\[
F_3(T) \equiv \int_0^T Z^3(t)\d t = O(T^c)?\
\leqno(3.5)
\]
Perhaps even $c = 3/4+\e$ is permissible? Is $c < 3/4$
impossible in} (3.5)?

\medskip
To this day the problem remains open.
However, if one considers the
cubic moment of $|Z(t)|$, then it is known that
$$
T(\log T)^{9/4} \;\ll\;  \int_1^T|Z(t)|^3\d t \;=\int_1^T|\zt|^3\d t
\;\ll\; T(\log T)^{9/4},\leqno(3.6)
$$
which establishes the true order of the integral in question. However, obtaining
an asymptotic formula for this integral remains a difficult problem.
The lower bound in (3.6) follows from general results of K. Ramachandra
(see his monograph \cite{Ram}), and the upper bound is a recent result of S. Bettin, V.
Chandee and M. Radziwi{\l}{\l} \cite{Bet2}.

\medskip
In \cite{Iv7}, equation (11.9), an explicit formula for the cubic
moment of $Z(t)$ was derived. This is
\[
\int\limits_T^{2T}Z^3(t)\d t= 2\pi\sqrt{\frac{2}{3}}
\sum_{(\frac{T}{2\pi})^{3/2}\le n\le (\frac{T}{\pi})^{3/2}}
d_3(n)n^{-\frac{1}{6}}\cos\bigl(3\pi n^{\frac{2}{3}}+\txt{\frac{1}{8}}\pi\bigr)
+O_\e(T^{3/4+\e}),\leqno(3.7)
\]
where as usual $d_3(n)$ is the divisor function
$$
d_3(n)\= \sum_{k\ell m=n}1\qquad(k,\ell,m,n\in \mathbb N),
$$
generated by $\z^3(s)$ for $\Re s >1$. Various techniques were used in \cite{Iv7}
to estimate the exponential sum in (3.7), but nothing better than the weak
$O_\e(T^{1+\e})$ seems to come out.

\medskip
A strong conjecture of the author is that
$$
\int_1^TZ^3(t)\d t \;=\;O_\e(T^{3/4+\e}).\leqno(3.8)
$$
Note that (3.8) would follow (by partial summation) from (3.7) and the bound
$$
\sum_{N<n\le N'\le2N}d_3(n)e^{3\pi in^{2/3}} \ll_\e N^{2/3+\e}.\leqno(3.9)
$$
It may be remarked that the exponential sum in (3.9) is ``pure'' in the sense that the function in
the exponential  does not depend on any parameter as, for example, the sum
$$
\sum_{N<n\le N'\le2N}n^{it}  \;=\;\sum_{N<n\le N'\le2N}e^{it\log n}
\qquad(1 \le N\ll\sqrt{t}\,),
$$
which appears in the approximation to  $\zt$ (see e.g., Theorem 4.1 of \cite{Iv1}),
depends on the parameter $t$.
However, the difficulty in the estimation of the sum in (3.9) lies in the presence of the
divisor function $d_3(n)$ which, in spite of its simple appearance, is quite difficult
to deal with.

\medskip
Finally we note that not much can be said about $F_{2k-1}(T)$ when $k\ge 3$. Even the conjecture
in (3.2) of Problem 1 remains open.
\medskip
\subsection{Moments of $Z(t)$ with shifts}

\medskip
A related and interesting problem  is to investigate
integrals of $Z(t)$ with ``shifts'', i.e., integrals where one (or more) factor $Z(t)$ is
replaced by $Z(t+U)$. The parameter $U$, which does not depend on the variable of
integration $t$, is supposed to be positive and $o(T)$ as $T\to\infty$, where $T$
is the order of the range of integration.

\medskip
Some results on such integrals already exist in the literature. For example,
R.R. Hall \cite{Hal} proved that, for $U = \a/\log T, \a\ll 1$,
we have uniformly
\begin{align*}
\int_0^TZ(t)Z(t+U)\d t &= \frac{\sin\a/2}{\a/2}T\log T + (2\gamma-1-2\pi)T\cos\a/2\\&
+ O\left(\frac{\a T}{\log T} + T^{1/2}\log T\right).\tag{3.10}
\end{align*}
Here $\gamma= -\G'(1) = 0.5772157\ldots\,$ is Euler's constant.
M. Jutila \cite{Jut3} obtained recently an asymptotic formula for the the  integral in (3.10) when
$0 < U \ll T^{1/2}$.

\medskip
S. Shimomura \cite{Shi} dealt with the quartic moment
$$
\int_0^T Z^2(t)Z^2(t+U)\d t,\leqno(3.11)
$$
under certain conditions on the real parameter $U$, such that $(|U| + 1)/\log T\to0$ as
$T\to\infty$. When $U \to 0+$, Shimomura's expression for (3.11) reduces to
$$
\int_0^T|\zt|^4\d t = \int_0^T Z^4(t)\d t = \frac{1}{2\pi^2}T\log^4T + O(T\log^3T).
\leqno(3.12)
$$
The (weak) asymptotic formula (3.12) is a classical result of A.E. Ingham \cite{Ing} of 1928.

\medskip
Finally  we mention that the author \cite{Iv9} obtained an asymptotic formula for the integral
of $Z^2(t)Z(t+U)$. This is formulated as

\medskip
{\bf Theorem 2}.
{\it For $\,0 < U = U(T) \le T^{1/2-\e}$ we have, uniformly in $U$},
\begin{align*}&
\int_{T/2}^{T}Z^2(t)Z(t+U)\d t = O_\e(T^{3/4+\e}) \,+
\\& \,+
2\pi\sqrt{\frac{2}{3}}\sum_{T_1\le n\le T_0}h(n,U)n^{-1/6+iU/3}
\exp(-3\pi in^{2/3}-\pi i/8)\Bigl\{1 + K(n,U)\Bigr\}.\tag{3.13}
\end{align*}
{\it Here} ($d(n)$ {\it is the number of divisors of $n$})
$$
h(n,U) := n^{-iU}\sum_{\delta|n}d(\delta)\delta^{iU},
\; T_0 := {\left(\frac{T}{2\pi}\right)}^{3/2}, \;T_1 := {\left(\frac{T}{\pi}\right)}^{3/2},
\leqno(3.14)
$$
$$
K(n,U) := d_2U^2n^{-2/3} + \cdots + d_kU^kn^{-2k/3}
+ O_k(U^{k+1}n^{-2(k+1)/3})\leqno(3.15)
$$
{\it for any given integer $k\ge 2$, with effectively computable constants
$d_2, d_3,\ldots\;$.}

\medskip
The interval of integration is $[T/2, \,T]$, since if it is $[0,T]$, then
$K(n,u)$ is not necessarily small. Note that, as $U\to0+$, the main term in (3.13)
becomes the main term in (3.7). In other words, Theorem 2 is a generalization of (3.7).
Therefore we may ask similar questions as was done in the case of $F_3(T)$.

\medskip
{\bf Problem 3}. {\it Is it true that there exists a constant $0 < c < 1$ such that the
integral in} (3.13) {\it is $O(T^c)$ uniformly for $0 < U = U(T) \le T^{1/2-\e}$}?

\medskip
The initial step in the proof of Theorem 2 is to write
$$
\int\limits_{T/2}^T Z^2(t)Z(t+U)\d t
= \frac{1}{i}\int\limits_{1/2+iT/2}^{1/2+iT}\z^2(s)\z(s+iU){(\chi^2(s)\chi(s+iU))}^{-1/2}\d s.
\leqno(3.16)
$$
The procedure of writing a real-valued integral like a complex integral is fairly
standard in analytic number theory. For example,  see
 the proof of Theorem 7.4 in E.C. Titchmarsh's monograph \cite{Tit}
on $\z(s)$ and M. Jutila's recent work \cite{Jut3}. It allows one flexibility
by suitably deforming the contour of integration
in the complex plane. Incidentally, this method of proof is different from the
proof of (3.7) in \cite{Iv7}, which is based on the use of approximate functional
equations.

In the  complex integral in (3.16) we replace the segment of integration by
$[1+\e + \hf iT, 1+\e + iT]$, and use the functional equation $\z(s) = \chi(s)\z(1-s)$
etc.
The problem is eventually reduced to the evaluation of exponential integrals
 whose saddle
point satisfies (when $U >0$) a non-trivial cubic equation
(i.e., $x^3 = ax+b$), whose solution is best found asymptotically. Lemma 1 is used
for the evaluation of the ensuing saddle points, and Theorem 2 follows eventually.

\medskip
\section{The distribution of values of $Z(t)$}

\medskip
Let $H = T^\t, 0 < \t\le1$, and
$$
I_+(T,H) := \int_{T,Z(t)>0}^{T+H}Z(t)\d t, \quad I_-(T,H) := \int_{T,Z(t)<0}^{T+H}Z(t)\d t.
\leqno(4.1)
$$
Also let
\begin{align*}
&{\cal J}_+(T,H) := \Bigl\{T< t\le T+H\,:\, Z(t)>0\Bigr\},\; {\cal K}_+(T,H) =
\mu\Bigl({\cal J}_+(T,H)\Bigr) \\&
{\cal J}_-(T,H) :=\Bigl\{T< t\le T+H\,:\, Z(t)<0\Bigr\}, \;{\cal K}_-(T,H) =
\mu\Bigl({\cal J}_-(T,H)\Bigr),
\tag{4.2}
\end{align*}
where $\mu(\cdot)$ denotes measure. We are interested in bounding $I_\pm, {\cal K}_\pm$.
In \cite{Iv5} the author proved that, unconditionally,
\begin{align*}
T(\log T)^{1/4} &\;\ll\; I_+(T,T)\ll T(\log T)^{1/4},\\
 T(\log T)^{1/4}\; &\ll\; -I_-(T,T)\ll T(\log T)^{1/4}.\tag{4.3}
\end{align*}

\medskip
{\bf Problem 4}. {\it Are there  constants $A_1, A_2>0$ such that}
$$
I_+(T,T) = (A_1 + o(1))T(\log T)^{1/4},  \; -I_-(T,T) = (A_2 + o(1))T(\log T)^{1/4}
\quad(T\to\infty)?
$$

 \medskip
 We present now a new result, which is contained in

 \medskip
 {\bf Theorem 3}. {\it Let $H = T^\t$ with $1/4 \le \t \le 1$. If the Riemann hypothesis
 is true, then for any number $k>1$ we have
 \begin{align*}
 {\cal K}_+(T,H) &\;\gg_k\; H(\log T)^{-k/4}\\
 {\cal K}_-(T,H) &\;\gg_k\; H(\log T)^{-k/4},
 \tag{4.4}
\end{align*}
where ${\cal K}_+(T,H), {\cal K}_-(T,H)$ are defined by} (4.2).

\medskip
{\bf Proof of Theorem 3}. First note that in \cite{Iv5}, \cite{Iv7}  the author,
for the left-hand sides in (4.4),
 obtained unconditionally the bounds $T(\log T)^{-1/2}$ when $H=T$.
The improvement in Theorem 3 is thus conditional, but the result holds
in a much more general case.

Assume the Riemann hypothesis (all complex zeros
of $\z(s)$ have real parts equal to 1/2). K. Ramachandra \cite{Ram} proved that
$$
\int_T^{T+H}|Z(t)|^k\d t \gg_k H(\log H)^{k^2/4}
\qquad(\log\log T\ll H \le T,\, k\in {\mathbb R}, k\ge0)
\leqno(4.5)
$$
%In fact, (4.5) for $k=1$ holds unconditionally.
An explicit value for the constant implicit in the $\gg$-symbol in (4.5) is to be found
in the work of M. Radziwi{\l}{\l} and K. Soundararajan \cite{RaSo}.
As for the upper bound for the integral in (4.5), we have
$$
\int_T^{T+H}|Z(t)|^k\d t \ll_k H(\log H)^{k^2/4}\qquad(T^\t\ll H \le T).\leqno(4.6)
$$
with $H = T^\t$ and $0<\t\le1$. This follows if one combines the results of A. Harper
\cite{Har} and the author \cite{Iv4}, both which are based on the method of
K. Soundararajan \cite{Sou}. Therefore it follows that
$$
H(\log T)^{1/4} \ll \int_T^{T+H}|Z(t)|\d t = \int_{{\cal J}_+(T,H)}Z(t)\d t -
\int_{{\cal J}_-(T,H)}Z(t)\d t.
$$
On the other hand,
$$
\int_{{\cal J}_+(T,H)}Z(t)\d t + \int_{{\cal J}_-(T,H)}Z(t)\d t
= \int_T^{T+H}Z(t)\d t \ll T^{1/4}
$$
by (3.4). Hence for $H = T^\t, 1/4 \le \t \le 1$ we have
\begin{align*}
H(\log T)^{1/4} &\ll \int_{{\cal J}_+(T,H)}Z(t)\d t,\\
H(\log T)^{1/4} &\ll -\int_{{\cal J}_-(T,H)}Z(t)\d t.
\end{align*}
Suppose now that $k$ satisfies $1<k<2$. By H\"older's inequality for integrals
and (4.6) we have
\begin{align*}
H(\log T)^{1/4} & \ll \int_{{\cal J}_+(T,H)}Z(t)\d t\le
\left(\int_{T,Z(t)>0}^{T+H}Z^{k}(t)\d t\right)^{1/k}\Bigl({\cal K}_+(T,H)\Bigr)^{1-1/k}
\\&
\ll_k (H(\log T)^{k^2/4})^{1/k}\Bigl({\cal K}_+(T,H)\Bigr)^{1-1/k}.
\end{align*}
This gives
\begin{align*}
H^{(k-1)/k}(\log T)^{(1-k)/4} &\ll_k \Bigl({\cal K}_+(T,H)\Bigr)^{(k-1)/k},\\
{\cal K}_+(T,H) &\gg_k H(\log T)^{-k/4}.
\end{align*}
In a similar fashion it is found that
\[
{\cal K}_-(T,H) \;\gg_k\; H(\log T)^{-k/4}.
\]
This completes the proof of Theorem 3.

\medskip
{\bf Problem 5}. {\it Do there exist positive constants $D_+, D_-$ such that
 $$
{\cal K}_+(T,T) = (D_+ + o(1))T,\quad  {\cal K}_-(T,T) = (D_- + o(1))T\quad(T\to\infty)?
 $$
 Is it true that} $D_+ = D_- = 1/2$?

 \medskip
 Of course, in general either ${\cal K}_+(T,H) \gg H$ or ${\cal K}_+(T,H) \gg H $ holds,
 but one cannot say which one of these lower bounds holds.

\medskip
 To continue our discussion on the evaluation of ${\cal K}_+(T,T)$
 (and ${\cal K}_-(T,T)$),  assume now the  Riemann Hypothesis
 and the simplicity of zeta zeros.
These very strong conjectures seem to be independent in the sense that it is not
known whether either of them implies the other one. Then (since $Z(0) = -1/2$) we have
$$
{\cal K}_+(T,T) = \mu\Bigl\{ T < t \le 2T : Z(t)>0\Bigr\} = \sum_{T<\gamma_{2n}\le 2T}(\gamma_{2n}
- \gamma_{2n-1}) + O(1),\leqno(4.7)
$$
where $0 < \gamma_1 < \gamma_2 < \ldots$ are the ordinates of complex zeros of $\z(s)$. Thus the
problem of the evaluation is reduced to the evaluation of the sum in (4.7). It seems reasonable
that the differences $\gamma_{2n} - \gamma_{2n-1}$ and  $\gamma_{2n+1} - \gamma_{2n}$
are evenly distributed, which heuristically
indicates that $D_+$ exists and that $D_+ = D_- = 1/2$. However, proving this is hard.

\medskip
Finally, to conclude our discussion on the distribution of values of
$Z(t)$, note that the  sum in (4.7) is related to the sum ($\a\ge0$ is fixed)
$$
\sum\nolimits_\a(T) : = \sum_{\gamma_n\le T}{(\gamma_n-\gamma_{n-1})}^\a,
$$
which was investigated in \cite{Iv21}. The sum $\sum_\a(T)$ in turn  can be connected to the
Gaussian Unitary Ensemble hypothesis (see A.M. Odlyzko \cite{Odl1}, \cite{Odl2})
and the pair correlation conjecture of H.L. Montgomery \cite{Mon}. Both of these conjectures
assume the Riemann Hypothesis and e.g., the former states that, for
$$
0 \le \a < \b < \infty, \quad \delta_n = \frac{1}{2\pi}(\gamma_{n+1}-\gamma_n)\log\left(
\frac{\gamma_n}{2\pi}\right),
$$
we have
$$
\sum_{\gamma_n\le T,\delta_n\in[\a,\b]}1 = \left(\int_\a^\b p(0,u)\d u + o(1)\right)
\frac{T}{2\pi}\log\left(\frac{T}{2\pi}\right)\qquad(T\to\infty).
$$
Here $p(0,u)$ is a certain probabilistic density, given by complicated functions
defined in terms of prolate spheroidal functions. In fact, in
\cite{Iv21} the author proved that,
if the RH and the Gaussian Unitary Ensemble hypothesis hold, then for $\a\ge0$ fixed
and $T\to\infty$,
$$
\sum\nolimits_\a(T) = \Bigl(\int_0^\infty p(0,u)u^\a\d u + o(1)\Bigr)
{\left(\frac{2\pi}{\log\bigl(\frac{T}{2\pi}\bigr)-1}\right)}^{\a-1}T.
$$
\smallskip
Also note that, since $\Re\log\zt
= \log|Z(t)|$, a classical result of A. Selberg (see \cite{Sel}, Vol. 1)
gives, for any real $\a<\b$,
$$
\lim_{T\to\infty}\frac{1}{T}\,\mu\Biggl\{\,t : t\in [T,2T],\,
\a < \frac{\log|Z(t)|}{\sqrt{\hf\log\log T}} < \b\,\Biggr\}
= \frac{1}{\sqrt{2\pi}}\int_\a^\b  e^{-\frac{1}{2}x^2}\d x,
$$
but here we are interested in the distribution of values of $Z(t)$ and not $|Z(t)|$.
\smallskip
Recently J. Kalpokas and J. Steuding \cite{KaSt} proved that for
$\phi\in [0,\pi)$,
$$
\sum_{0<t\le T,\,\zeta(\frac{1}{2}+it)\,\in\, e^{i\phi}{\mathbb R}}\zt = \Bigl(2e^{i\phi}
\cos\phi\Bigr) \frac{T}{2\pi}\log\frac{T}{2\pi e} + O_\e(T^{1/2+\e}),\leqno(4.8)
$$
and an analogous result holds for the sums of $|\zt|^2$. It is unclear whether (4.8)
and the other approaches mentioned above can be put to
use in connection with our problems.

 \medskip
 \vfill
 \eject

\end{document}